\documentclass[a4paper,11pt,leqno]{amsart}

\usepackage{amssymb,hyperref}
\usepackage[latin1]{inputenc}
\usepackage[francais,english]{babel}
\usepackage{graphicx}
\usepackage{color}

\usepackage[T1]{fontenc}
\usepackage[latin1]{inputenc}
\usepackage[all]{xy}
\usepackage{stmaryrd}
\usepackage{epsfig,psfig}
\usepackage{float}

\usepackage{psfrag}
\usepackage{wrapfig}
\usepackage{amsmath}
\usepackage{amssymb,latexsym}
\usepackage[mathscr]{eucal}
\usepackage{url}
\usepackage{fancybox}
\usepackage{bbm} 

\swapnumbers

\numberwithin{equation}{section}

\DeclareMathAlphabet{\mathrmsl}{OT1}{cmr}{m}{sl}

\newcommand{\oper}[3][n]{\newcommand{#2}{\mathop{\mathrm{#3}}%
\ifx n#1\nolimits\else\limits\fi} }
\newcommand{\rsoper}[3][n]{\newcommand{#2}{\mathop{\mathrmsl{#3}}%
\ifx n#1\nolimits\else\limits\fi} }

\newcounter{mnotecount}[section]
\renewcommand{\themnotecount}{\thesection.\arabic{mnotecount}}
\newcommand{\mnote}[1]
{\protect{\stepcounter{mnotecount}}$^{\mbox{\footnotesize  $
      \bullet$\themnotecount}}$ \marginpar{\raggedright\tiny\em
    $\!\!\!\!\!\!\,\bullet$\themnotecount: #1} }

\newcommand{\eq}[1]{\ref{#1}}

\renewcommand{\qed}{\hfill $\blacksquare$ \medskip \\}

\newcommand{\question}{\noindent {\bf Question:\ }}

\theoremstyle{main}
\newtheorem{thm} {\bf  Theorem} [section]
\newtheorem{coro} [thm] {\bf  Corollary}
\newtheorem{lem} [thm] {\bf  Lemma}

\newtheorem{d\'ef}[thm]{\bf Definition}
\newtheorem{Rq}[thm]{\bf  Remark }
\newtheorem{Rqs}[thm]{\bf  Remarks}
\newtheorem{Exemple} [thm] {\bf  Example}

\DeclareFontFamily{OT1}{rsfs}{}
\DeclareFontShape{OT1}{rsfs}{m}{n}{ <-7> rsfs5 <7-10> rsfs7 <10-> rsfs10}{}
\DeclareMathAlphabet{\mycal}{OT1}{rsfs}{m}{n}

\newcommand{\dVol}{\operatorname{dVol}}
\newcommand{\Vol}{\operatorname{Vol}}

\newcommand{\Ric}{\operatorname{Ric}}

\newcommand{\Supp}{\operatorname{Supp}}
\newcommand{\R}{{\Bbb R}}

\newcommand{\N}{{\Bbb N}}

\renewcommand{\l}{\lambda}
\newcommand{\x}{  \textbf{\textit{x}}}

\title[eigenvalues estimates for Neumann]{Eigenvalues estimate for the Neumann problem of a bounded domain}
\author{Bruno COLBOIS and Daniel MAERTEN}
\date{\today}
\address{Université de Neuchâtel, Institut de Math\'ematiques, rue Emile Argand 11,  CH--2007 Neuch\^atel, Suisse.}
\address{Université de Tours, LMPT, Fédération Denis Poisson, Faculté des Sciences, Parc de Grandmont, F--37200 Tours, France.}
\email{{\color{red}bruno.colbois@unine.ch\\
daniel.maerten@yahoo.fr}}

\begin{document}

\maketitle

\begin{abstract}
In this note, we investigate upper bounds of the Neumann eigenvalue problem for the Laplacian 
of a domain $\Omega$ in a given complete (not compact a priori) Riemannian manifold $(M,g)$.  
For this, we use test functions for the Rayleigh quotient subordinated to a family of open sets
constructed in a general metric way, interesting for itself. As applications,
we prove that if the Ricci curvature of $(M,g)$ is bounded below $\Ric^{g}\ge -(n-1)a^2$, $a \ge 0$, 
then there exist constants $A_n > 0, B_n>0$ only depending on the dimension, such that  
$$	\l_k(\Omega)\le A_n a^2 + B_n\left(\frac{k}{V}\right)^{2/n}\ ,$$
where $\l_k(\Omega)$ $(k\in\N^{*})$ denotes the $k$--th eigenvalue of the Neumann problem on any bounded domain $\Omega\subset M$ of volume ${V=\Vol(\Omega,g)}$. Furthermore,  this upper bound is clearly in agreement with the Weyl law. As a corollary, we get also an estimate which is analogous to Buser's upper bounds of the spectrum of a compact Riemannian manifold with lower bound on the Ricci curvature.
\end{abstract}

\section{Introduction}

The goal of this paper is to give upper bounds for the spectrum of the Laplacian acting 
on compact domains of given volume of a complete Riemannian manifold with Ricci curvature bounded below, and, as far as possible, 
to make these estimates optimal with respect to the Weyl law. 

For compact Riemannian manifolds without boundary, the following result was 
proved by P. Buser in \cite{Bu1} (Satz 7), \cite{Bu2} (Thm. 6.2 (c)) (see also Li-Yau in \cite{PL1} (Thm.16)). If
$\{\lambda_k\}_{k=1}^{\infty}$ denote the spectrum of the Laplacian acting on functions, then:

\begin{thm}\label{buser} 
Let $(M^{n},g)$ be a compact $n$-dimensional Riemannian manifold with Ricci curvature bounded below $\Ric^{g} \ge -(n-1)a^2$, $a \ge 0,$ 
and of volume $V$. 
 
There exists a constant $C_n\ge 1$ only depending on the dimension, such that for all $k\in \N^{*}$, we have  
\begin{equation}\label{coroneg}
	\l_k(M,g)\le \frac{(n-1)^2}{4} a^2+ C_n\left(\frac{k}{V}\right)^{2/n}\ .
\end{equation}
 
\end{thm}

\begin{Rqs}
\begin{enumerate}
\item
In \cite{PL1}, the constant $C_n$ depends also on the diameter.

\item
In dimension higher than 2, a normalization on the volume is not enough to control the spectrum: namely, on any compact manifold of dimension higher than 2,
one can find a metric of given volume,  with arbitrarily large first non--zero eigenvalue $\l_{2}$ of the Laplacian, 
in vertue of the result of B.~Colbois and  J.~Dodziuk \cite{CD}.  

\item
When $\Ric^{g} \ge 0$, we deduce that there exists $C_n>1$ with $\l_k(M,g)\le C_n\left(\frac{k}{V}\right)^{2/n}$  for all $k$. 
However, when Ricci is not
supposed positive, then the presence of a term like $\frac{(n-1)^2}{4} a^2$ is necessary: by a result of R.~Brooks~\cite{Br}, 
it is possible to find a family of compact hyperbolic manifolds with volume going to infinity and a positive uniform lower bound on 
the first nonzero eigenvalue.  
\end{enumerate}
\end{Rqs}
 
The idea of the proof of Theorem \ref{buser} is to consider $k$ disjoint balls of radius $r$
which almost cover the manifold $(M,g)$, with $r$ around $\left(\frac{V}{k}\right)^{1/n}$, 
and to apply then Cheng's theorem \cite{Ch}. However, such a theorem does not exist on manifolds with boundary, 
and with Neumann boundary condition.
A reason for this is that there is no Bishop-Gromov theorem: indeed, even for a Euclidean domain, it is not possible to
control the volume of a ball of radius $2r$ with respect to the volume of a ball of radius $r$ and same center. See
also Example 1.4 in \cite{Bu2}.

\smallskip
This does not mean that a result in the spirit of Theorem \ref{buser} does not exist for domains.
Namely, P.~Kröger \cite{K} proved thanks to harmonic analysis, that on bounded 
Euclidean domains, the $k$--th eigenvalue of the Neumann problem was bounded by above by some expression $
C_{n}\left(k/\left|\Omega\right|\right)^{n/2},$ where $C_{n}$ only depends upon the dimension. 
An analogous result can be derived from the much more general and difficult work of N.~Korevaar \cite{Korevaar} (
see also \cite{GY}), for bounded domains of non--negative Ricci curvature  
manifolds, and also for bounded domains of negative Ricci curvature compact 
manifolds (in this case the bound depends on the diameter).

\medskip

This naturally leads to the\\

\question 
What can be said for bounded domains of a complete Riemannian manifold with Ricci curvature bounded below ?\\

In this note, we consider the Neumann eigenvalue problem for the Laplacian of a bounded domain $\Omega$ with smooth boundary, in a given complete (not compact a priori) Riemannian manifold $(M,g)$. 
More precisely, we search for a couple $(\lambda,u)\in \R \times C^{\infty}\left(\overline{\Omega}\right)$ which is a solution of the following boundary elliptic  problem 
$$\left\{
\begin{array}{ll}
	\Delta u =\l u &  \textrm{ on } \Omega \\
	\frac{\partial u }{\partial \nu}=0 &  \textrm{ on } \partial \Omega \ ,
\end{array}\right. $$
where $\Delta$ is the non--negative Laplacian of the metric $g$ and $\nu$ the outward unit normal of $\partial \Omega$. Since $\Omega$ is bounded with smooth boundary, the spectrum of $\Delta$ on $\Omega$ is an unbounded sequence of real numbers $\left(\l_{k}(\Omega)\right)_{k\in\N^{*}}$ which can be increasingly ordered 
$$ 0= \l_{1}(\Omega) < \l_{2}(\Omega) \le \cdots \le  \l_{k}(\Omega)\le \l_{k+1} (\Omega)\le \cdots  \quad .$$
There exist standard variational characterisations of the spectrum of $\Delta$ which can be found for instance in the book of  P.~Bérard~\cite{Be} (or in \cite{GHL}).\\

The main result of this article is the following.
 
\begin{thm} \label{t1} 
Let $(M^{n},g)$ be a complete $n$-dimensional Riemannian manifold with Ricci curvature bounded below $\Ric^{g} \ge -(n-1)a^2$, $a \ge 0$. 

\smallskip
There exist constants $A_n > 0, B_n>0$ only depending on the dimension, such that for all $k\in \N^{*}$, $V>0$ 
and for each bounded domain $\Omega \subset M$, with smooth boundary and volume $V$, we have  
\begin{equation}\label{riccineg}
	\l_k(\Omega)\le A_n a^2 + B_n\left(\frac{k}{V}\right)^{2/n}\ .
\end{equation}

\end{thm}

If the manifold $M$ is compact, an interesting special case is to choose $\Omega =M$, and we recover Theorem~\ref{buser}, up to the value of the constant $A_n$ which is not
equal to $\frac{(n-1)^2}{4}$ in our paper.

\medskip
The proof Theorem~\ref{t1} goes in the same spirit as the proof of Theorem \ref{buser}: in order to bound $\lambda_k(\Omega)$,
we consider $k$ disjoint sets $A_1,...,A_k$ in $\Omega$ of measure of the order of $\frac{Vol(\Omega)}{k}$,
 and introduce test functions $f_1,...,f_k$ subordinated to these sets. 
We estimate the Rayleigh quotient of these fonctions by a direct calculation, which gives the theorem. The main improvement
of this paper is the construction of an adapted family of sets $A_1,..,A_k$, more convenient for our purpose as balls.
As this construction is interesting by itself and will be used in other contexts, we present it in a rather abstract (indeed metric) way.

\medskip
The paper is organised as follows: the metric construction of our sets is done in Section~\ref{s2}, 
and in Section~\ref{s3} we will use them so as to prove Theorem~\ref{t1} by producing some test functions 
for the variational characterisation of the spectrum.

\section{A metric approach} \label{s2} 

In this section, we formalize the geometric situation of Theorem~\ref{t1} (a bounded domain in a complete manifold) in a more general setting (a bounded domain in a complete metric space). 
More precisely, let  $(X,d)$ be a complete, locally compact metric space,  ${Y \subset X}$ a bounded Borelian subset endowed with the induced distance, 
and $\mu$ a Borelian measure with support in $\overline{Y}$ such that $\mu(Y)=\omega$,  $ 0 < \omega < \infty$. We will need in addition the following 
technical assumptions: 
\begin{itemize}
	\item[(H1)] For each $r>0$, there exists a constant $C(r)>0$ such that each ball of radius $4r$ in $X$ may be covered by $C(r)$ balls of radius $r$. Moreover, $r\mapsto C(r)$ is an increasing function of the radius.\\
	\item[(H2)] We suppose that the volume of the $r$--balls tends to 0 uniformly on $X$, namely $\underset{r\rightarrow0}{\lim } \sup\{\mu(B(x,r)):x\in X\}=0$.  
  However, taking (H1) into account, this volume condition is equivalent to ${\underset{r\rightarrow0}{\lim } \sup\{2C(r)\mu(B(x,r)):x\in X\}=0}$ which is the (more convenient) condition that will be used in the remainder of the article.
\end{itemize}
It is important to remark that these hypothesis are quite natural since they make part of the metric properties of the Riemannian manifolds that are involved in Theorem~\ref{t1}. 
These specific metric properties are collected in the following fundamental example.

\begin{Exemple}\label{exemple} A typical example of a couple $(X,Y)$ satisfying the hypothesis {\rm (H1),(H2)}  
is to choose  $X$ as a complete $n$-dimensional Riemannian manifold $(M,g)$ with Ricci curvature 
bounded below $\Ric^{g} \ge -(n-1)a^2$, $a \ge 0$ (which are the class of manifolds involved in Theorem~\ref{t1}), and as $Y$ a bounded domain with smooth boundary in $M$. The distance $d$ is the distance associated to the Riemannian metric $g$, the measure $\mu$ is the restriction  to $Y$ of the Riemannian measure of $g$. The existence of the constant $C(a,r)$ is given by the classical  Bishop-Gromov inequality thanks to the lower bound on the Ricci curvature of $g$ (see \cite{Sa} p. 156). Precisely, for $0 <r<R$, and for each point $p\in M$, we have 
\begin{equation}\label{bishop}
\frac{\Vol(B(p,R),g)}{\Vol(B(p,r),g)} \le \frac{v_a(R)}{v_a(r)} \ , 
\end{equation}
where $v_a(R)$ denotes the volume of a ball of radius $R$ in $\Bbb M^{n}_{a}$, 
the simply connected $n$--dimensional manifold of constant sectional curvature $-a^2$.

This gives a bound on the number of  balls of radius $r$ that are necessary to cover a ball of radius $4r$ 
(this property known as the paking lemma is a consequence of Inequality~(\ref{bishop})). 
In fact, fix  $B_{4r}$ a $4r$--ball and consider $\left\{B(x_{i},r/2)\right\}_{i\in I}$ a maximal family of disjoint balls whose center $x_{i}$ live in $B_{4r}$; then the corresponding family of $r$--balls $\left\{B(x_{i},r)\right\}_{i\in I}$ cover $B_{4r}$. In consequence, we can cover a ball of radius $4r$ with
$\le 1 + \left[\frac{v_{a}(4r+r/2)}{v_{a}(r/2)}\right]$ $r$--balls. We just define 
$$C(a,r)=\max_{t \le r}\left\{1 + \left[\frac{v_{a}(4t+t/2)}{v_{a}(t/2)}\right]\right\} \ .$$ 
The increasing character of $r\mapsto C(a,r)$ is by definition. 
 
\medskip
Furthermore, as $r \longrightarrow 0$, the ratio $\frac{\Vol(B(p,r),g)}{v_a(r)} \longrightarrow 1$, we obtain
$$\Vol((B(p,R),g)\le v_a(R) \ ,$$
 and consequently $\mu(B(p,r)):=\Vol(B(p,r)\cap Y,g)$ goes uniformly to $0$ as $r\to 0$.
\end{Exemple}

We prove in the sequel that, under our technical assumptions, one can build some subsets $A$ and $D$  satisfying certain volume conditions.

\begin{lem}\label{l1} 
Let $(X,d)$ be a complete, locally compact metric space,  ${Y \subset X}$ a bounded Borelian with the induced distance, and $\mu$ a Borelian measure 
with support in $\overline{Y}$ such that $\mu(Y)=\omega$,  $ 0 < \omega < \infty$ and $\mu(\overline{Y}\setminus Y)=0$.  In addition, we make the  hypothesis {\rm (H1),(H2)}.\\
Let $0< \alpha \le \frac{\omega}{2}$. Thanks to {\rm (H2)} there exists $r>0$ with $\sup\{2C(r)\mu(B(x,r)):x\in X\}\le \alpha$.\\
Then there exist $A,D \subset Y$ such that $A\subset D$ and 
$$\left\{
\begin{array}{l}
	\mu(A)\ge \alpha\\
	\mu(D) \le 2C(r)\alpha\\
	d(A, Y\cap D^{c}) \ge3r 
\end{array}\right. \ .$$
\end{lem}

\bigskip
\noindent
\textbf{Proof.} 
We fix the positive  numbers $r$ and $\alpha$. Let us consider any positive integer $m\in\N^{*}$ and define a non--negative application ${\Psi_{m}:X^{m}= \underbrace{X\times X \times \cdots \times X }_{m \  {\rm times}} \longrightarrow   \R }$ by the relation
$$\Psi_{m}:\x=\left(x^{j}\right)^{m}_{j=1} \longmapsto  \mu\left(\bigcup^{m}_{j=1} B\left(x^{j},r\right) \right) \ ,$$
which is simply the restriction of the measure $\mu$ to $\mycal{U}_{m}(r)$  a particular class of open sets which is defined by
$$\mycal{U}_{m}(r):=\left\{ \bigcup^{m}_{j=1} B\left(x^{j},r\right) /  \  \left(x^{j}\right)^{m}_{j=1}\in X^{m} \right\} \ .$$ 
Since $(X,d)$ is a complete and locally compact metric space, it is also the case of the finite product $X^{m}$ when it is endowed with the product distance. Then for each $m\in\N^{*}$ there exists some $\x_{\max, m} \in X^{m}$ (not necessary unique) such that
$$ \Psi_{m}(\x_{\max,m})= \max_{X^{m}} \Psi_{m} = \max_{\mycal{U}_{m}(r)}  \mu = \mu  \left(\bigcup^{m}_{j=1} B\left(x^{j}_{\max,m},r\right)\right) \ .$$
We first prove that there exists a finite integer  $k\in\N^{*}$ such that $\Psi_{k}(\x_{\max,k})\ge \alpha$ and  $\Psi_{k-1}(\x_{\max,k-1})\le \alpha$. Indeed, consider the function $\xi: \N^{*} \longrightarrow \R $ defined by the relation  ${\xi(m)=\Psi_{m}(\x_{\max,m})}$. On one hand, the condition $\sup\{2C(r)\mu(B(x,r)):x\in X\}\le \alpha$ obviously implies $\xi(1)\le \frac{\alpha}{2C(r)}\le \alpha $. On the other hand, since $\Supp \mu\subset \overline{Y}$,  there exists a radius $R>0$ large enough such that $\mu ( B(z,R))\ge 3\omega/4 $, for a certain $z\in X$. But it can be clearly deduced from Assumption (H1) that $B(z,R) $ can be finitely covered by  $m_{0}\in\N ^{*}$ balls of radius $r$ (notice that $m_{0}$ depends on $R$).  Consequently it turns out
$$	\frac{3\alpha}{2}  \le  \frac{3\omega}{4}\le \mu \left( B(z,R)\right)\le  \max_{\mycal{U}_{m_{0}}(r)} \Psi_{m_{0}}=  \xi(m_{0}) \ .$$
Thereby the function  $\xi: \N^{*} \longrightarrow \R $ satisfies $\xi(1)\le \alpha$ and $\xi(m_{0})  \ge\frac{ 3\alpha}{2}$, which entails the existence of  some  $k\in\N^{*}$ such that $\Psi_{k}(\x_{\max,k})\ge \alpha$ and  $\Psi_{k-1}(\x_{\max,k-1})\le \alpha$.\\
We now set $U_{k}:=\underset{1\le j \le k}{\bigcup} B\left(x^{j}_{\max,k},r\right) $ and $V_{k}:=\underset{1\le j \le k}{\bigcup} B\left(x^{j}_{\max,k},4r\right) $. The next step is to show that
$$ \mu(V_{k})\le C(r)\mu(U_{k}) \ .$$
Still according to Assumption (H1), $V_{k}$ is covered by $kC(r)$ balls of radius $r$,  namely $V_{k}\subset\underset{1\le j \le kC(r)}{\bigcup}B_{j}$, where the $B_{j}$ are balls of radius $r$. But it is quite clear that this union of $r$--balls can be written as $\underset{1\le j \le kC(r)}{\bigcup}B_{j}= \underset{1\le j \le C(r)}{\bigcup}W_{j}$ where each $W_{j}\in \mycal{U}_{k}(r)$. It follows
\begin{eqnarray*}
	\mu(V_{k})\le \mu\left( \underset{1\le j \le kC(r)}{\bigcup}B_{j} \right)& =& \mu\left(\underset{1\le j \le C(r)}{\bigcup}W_{j}\right) \\
	&\le & \sum^{C(r)}_{j=1}\mu(W_{j}) \\
	&\le & C(r) \max_{\mycal{U}_{k}(r)} \mu = C(r) \xi (k)= C(r) \mu(U_{k} ) \ .
\end{eqnarray*}
We finally define the  sets $A:= Y \cap U_{k}$ and $D:=Y \cap V_{k}$. We only have to check that they satisfy the properties stated in Lemma~\ref{l1}. We observe that $\mu(A)=\mu(U_{k})$ since the measure $\mu$ is supported in $\overline{Y}$ and $\mu(\overline{Y}\setminus Y)=0$. Besides, $U_{k}$ can be written as the union of an element of $\mycal{U}_{k-1}(r)$ and an element of $\mycal{U}_{1}(r)$ so that 
$$\mu(A)\le \xi(k-1)+ \xi (1)\le \alpha \left(1+ \frac{1}{2}\right) \ .$$
Still since $\Supp \mu =\overline{Y}$, we obtain $\mu(D)=\mu(V_{k})\le C(r)\mu(U_{k}) = C(r) \mu(A)\le 2 C(r) \alpha$. By the definition of $U_{k}$ and $V_{k}$, we straightforwardly have $d(A,Y\cap D^{c})\ge3r$. \qed

In section \ref{s3}, we will use the following corollary of Lemma \ref{l1} to make the proof of Theorem \ref{t1}. 
We give therein an explicite construction of the domains that were mentioned at the end of the introduction.

\begin{coro}\label{c1} 
Let $(X,d)$ be a complete, locally compact metric space,  ${Y \subset X}$ a bounded Borelian with 
the induced distance, and $\mu$ a Borelian measure with support in $\overline{Y}$ 
such that $\mu(Y)=\omega$,  $ 0 < \omega < \infty$ and $\mu(\overline{Y}\setminus Y)=0$.  In addition, we make the  hypothesis {\rm (H1),(H2)} as in Lemma~\ref{l1}, and take $N$ a positive integer.\\
Let $r>0$ such that $4C^2(r)\mu(B(x,r)) \le \frac{\omega}{N}$ holds for all $x \in X$, and let $\alpha= \frac{\omega}{2C(r)N}$.
Then, there exist $N$ measurable subsets $A_1,...,A_N \subset Y$ such that $\mu(A_i)\ge \alpha$ and, for each $i \not =j$, $d(A_i,A_j)\ge 3r$.
\end{coro}


\noindent
\textbf{Proof.} 
We construct the family $\left(A_{j}\right)^{N}_{j=1}$ by finite induction applying Lemma~\ref{l1}. 

\begin{itemize}
	\item[$\bullet\ j=1$.] 
We set $(X_{1},d_{1},\mu_{1})=(X,d,\mu)$ and $Y_{1}=Y$, which satisfy the assumptions of Lemma~\ref{l1}. Therefore there exist $A_{1},D_{1}$ 
such that $A_{1}\subset D_{1}\subset Y_{1}=Y$ and 
$$\left\{
\begin{array}{rll}
	\mu(A_{1})&  \ge & \alpha\\
	 \mu(D_{1}) & \le & 2C(r)\alpha = \frac{\omega}{N}\\
	d(A_{1},  Y_1\cap D_{1}^{c}) & \ge & 3r 
\end{array}\right. \ .$$\\
	\item[$\bullet\ j=2$.] 
We set $(X_{2},d_{2},\mu_{2})=(X,d,\mu_{|Y_{2}})$ and $Y_{2}=D^{c}_{1}\cap Y_1$, which satisfy the assumptions of Lemma~\ref{l1} with $\omega_{2}= \mu_{2}(Y_{2})\ge \omega\left(1-\frac{1}{N}\right) = \omega \left(\frac{N+1-2}{N}\right)\ge\alpha $. Therefore there exist $A_{2},D_{2}$ such that $A_{2}\subset D_{2}\subset Y_{2}=D^{c}_{1}\cap Y_1$ and 
$$\left\{
\begin{array}{rll}
	\mu(A_{2})&  \ge& \alpha\\
	\mu(D_{2}) & \le & 2C(r)\alpha = \frac{\omega}{N}\\
	d(A_{2}, Y_{2}\cap D_{2}^{c}) & \ge & 3r 
\end{array}\right. \ .$$
As $A_{1}\subset D_{1}$ and $A_{2}\subset Y_1 \cap D_{1}^{c}$ we get $d(A_{1},A_{2})\ge d(A_{1},Y_1 \cap D^{c}_{1}) \ge 3r$ thanks to the case $j=1$.\\
	\item[$\bullet  j\ge3$.] 
We suppose that we have already constructed the families
 $\left(A_{s}\right)^{j-1}_{s=1}$ and $\left(D_{s}\right)^{j-1}_{s=1}$ that satisfy the conditions
$$\left\{
\begin{array}{l}
	A_{s}\subset D_{s} \subset Y \cap \left( D_{1}\cup \cdots \cup D_{s-1}\right)^{c} = Y_s, \quad s\le j-1\\
	d(A_{s},A_{t}) \ge 3r \quad s\neq t ,\\
	\mu\left( D_{1}\cup \cdots \cup D_{j-1}\right) \le \omega \left(\frac{j-1}{N}\right) \ .
\end{array}\right. $$
We set $(X_{j},d_{j},\mu_{j})=(X,d,\mu_{|Y_{j}})$ and $Y_{j}=Y \cap \left( D_{1}\cup \cdots \cup D_{j-1}\right)^{c}$, which satisfy the assumptions of Lemma~\ref{l1} with ${\omega_{j}=\mu_{j}(Y_{j})\ge \omega\left(1-\frac{j-1}{N}\right) = \omega \left(\frac{N+1-j}{N}\right)\ge\alpha }$ if $j\le N$. Therefore there exist $A_{j},D_{j}$ such that $A_{j}\subset D_{j}\subset Y_{j}$ and 
$$\left\{
\begin{array}{rll}
	\mu(A_{j})&  \ge & \alpha\\
	\mu(D_{j}) & \le & 2C(r)\alpha  = \frac{\omega}{N}\\
	d(A_{j}, Y_{j}\cap D_{j}^{c}) & \ge & 3r 
\end{array}\right. \ .$$
As $A_{j}\subset Y \cap \left( D_{1}\cup \cdots \cup D_{j-1}\right)^{c}\subset Y\cap \left( D_{1}\cup \cdots \cup D_{s-1}\right)^{c}=Y_{s} $, $s<j$, and $A_{s}\subset D_{s}$,  we get $d(A_{j},A_{s})\ge d(A_{s},Y_{s}\cap D^{c}_{s}) \ge 3r$ thanks to the case $j=s$. As already said, we can proceed this construction so longer we have enough volume to do it, that is $N$ times.\qed
\end{itemize}

\section{Proof of Theorem~\ref{t1}.}\label{s3}
Let $(M^{n},g)$ be a complete $n$-dimensional Riemannian manifold with Ricci curvature bounded below $\Ric^{g}\ge -(n-1)a^2$, 
and $\Omega \subset M$ a bounded domain of volume $V$, with smooth boundary.\\

We observe first that, by renormalisation, it is enough to prove the theorem for the case $a=1$: namely, 
if Theorem~\ref{t1} is true for ${a=1}$, and if $g$ is a Riemannian metric with
${\Ric^g \ge -(n-1)t^2 g}$, then $g_0=t^2 g$ satisfies $\Ric^{g_0} \ge -(n-1) g_0$. 
Since we have $\lambda_k(g_0) \le A_n+B_n \left(\frac{k}{V(g_0)}\right)^{2/n}$, then, because $\lambda_k(g)=t^2\lambda_k(g_0)$ and
$V(g)=t^nV(g_0)$, we get $\lambda_k(g)\le A_nt^2+B_n\left(\frac{k}{V}\right)^{2/n}$.

\medskip
So, let use prove Theorem \ref{t1} for $a=1$.
As in Example~\ref{exemple}, let us consider the Borelian measure $\mu$ which is the restriction to 
the domain $\Omega$ of the Riemannian volume of $(M,g)$. 

\smallskip
In order to prove Theorem~\ref{t1}, we will use the classical variational characterization of the spectrum: 
to estimate $\l_k$ from above, 
it suffices to construct an $H^{1}(\Omega)$-orthogonal family of $k$ test functions 
$\left(f_{j}\right)^{k}_{j=1}$, such as each $f_{j}$ has controled Rayleigh quotient. 
In the sequel, we construct test functions with disjoint support related
to the sets $A_1,...,A_k$ arising from Corollary \ref{c1}, so that it immediately implies orthogonality in $H^{1}(\Omega)$.\\

\begin{lem}\label{l2} 
Let $A \subset M$ a subset as in Corollary \ref{c1}. Let $A^r:=\{x\in M: d(x,A) \le r \}$, $r>0$. There exists a function $f$ supported in $A^r$ whose restriction to $\Omega$ is of Rayleigh quotient 
$$R(f) \le \frac{1}{r^2}\frac{\mu(A^r\setminus A)}{\mu(A)}.$$
\end{lem}

\smallskip
\noindent\textbf{Proof.} Let us define a \emph{plateau} function
$$f(p)=\left\{
\begin{array}{cll}
	1 & {\rm if } & p\in A \\
	1-\frac{d(p,A)}{r} & {\rm if } & p\in \left(A^r\setminus A\right) \\
	0 & {\rm if } & p\in (A^r)^{c}  \ .
\end{array}\right.$$
In Corollary \ref{c1}, the domain $A$ is a finite union of metric balls and intersection with complement of balls. The boundary is not smooth, but the function $d(\partial A, \cdot )$ "distance to the boundary of $A$"  is well known to be 1--Lipschitz on $M$. According to Rademacher's theorem (see Section~3.1.2, page 81--84 in \cite{EG}), $d(\partial A, \cdot )$ is differentiable $\mycal{L}^{n}$  almost everywhere (since  $\dVol_{g}$ is absolutely continuous with respect to Lebesgue's measure $\mycal{L}^{n}$), and its $g$--gradient satisfies  $\left|\nabla d(\partial A, \cdot )\right|_{g}\le 1$, $\mycal{L}^{n}$  almost everywhere. It comes out that the gradient of $f$ satisfies $\mycal{L}^{n}$  almost everywhere
$$\left|\nabla f(p)\right|_{g}\le\left\{
\begin{array}{cll}
	\frac{1}{r} & {\rm if } & p\in \left(A^r\setminus A\right) \\
	0 & {\rm if } & p\in \left(A^r\setminus A\right)^{c}  \ .
\end{array}\right.$$
We immediately deduce 
$$R(f) = \frac{\int_{\Omega}\left|\nabla f\right|^{2}_{g}\dVol_{g}}{\int_{\Omega}f^{2}\dVol_{g}} \le \frac{1}{r^2}\frac{\mu(A^r\setminus A)}{\mu(A)}.$$
\qed

\bigskip
\noindent
\textbf{Proof of Theorem~\ref{t1}.}
As already said, we apply Corollary~\ref{c1}: let $k \in \N^{*}$ and set ${N=2k}$. 
As the volume of the $r$--balls uniformly tends to 0 (see assumption (H2)),  there exist
$r>0$  with $r$ small enough so that 
\begin{equation}\label{techniq}
	{2C(r)\mu(B(x,r)) \le \alpha:= \frac{V}{4C(r)k}} \ ,
\end{equation}
holds for every $x\in M$. Corollary~\ref{c1} gives the existence of $2k$ measurable subsets $A_1,...A_{2k}$ of measure 
$\mu(A_i) \ge \frac{V}{4C(r)k}$ with $d(A_i,A_j)\ge 3r$ if $i\not =j$. 
In particular, the corresponding sets $A_i^r$ and $A_j^r$ are also disjoint.\\
We can now apply the construction of Lemma~\ref{l2} and we get an $H^1(\Omega)$-orthogonal 
family of $2k$ test functions $\left(f_{j}\right)^{2k}_{j=1}$, of disjoint supports and whose Rayleigh quotient satisfies
$$R(f_i) \le \frac{1}{r^2}\frac{\mu(A_i^r\setminus A_i)}{\mu(A_i)}.$$
At this point, Corollary~\ref{c1} does not give any control on $\mu(A_i^r)$. Let 
$$Q= \sharp \left\{i\in \left\{1,...,2k \right\}: \mu(A_i^r) \ge \frac{V}{k}\right\}.$$
As $\Vol(\Omega,g)=V$, we already see that $Q \le k$, so that for at least $k$ of these $2k$ subsets $A_1,...,A_{2k}$, we have $\mu(A_i^r) \le \frac{V}{k}$.
We choose the corresponding functions as test functions. For such a function $f$, we have, 
as $\mu(A_i^r\setminus A_i)\le \frac{V}{k}$ and $\mu(A_i)\ge \alpha=\frac{V}{4C(r)k}$, that  
\begin{equation}\label{bornepreuve}
	R(f) \le \frac{1}{r^2} \frac{V/k}{V/4C(r)k}=\frac{4C(r)}{r^2} \ .
\end{equation}

Our aim is now to prove an upper bound of the kind
$$\lambda_k(g) \le A_n + B_n \left(\frac{k}{V}\right)^{2/n} \ .$$
Let $\omega'_n>0$  the positive 
constant such that $\mu\left(B(x,r)\right)\le \omega'_n r^{n}$ for radius $r\le1$ in the hyperbolic space of curvature $-1$. 
We then define the integer $k_0 = \left[\frac{V}{8C(1)^2 \omega'_n }\right] +1 $ (remark that it strongly depends on the volume) and for every $k\ge k_0$, we set
   	$$ r_k = \left(\frac{V}{k}  \frac{1}{8 C(1)^{2} \omega'_n }\right)^{1/n} \ .$$
Clearly, $r_k\le 1$ and (\ref{techniq}) holds, since by definition ${8C(r_k)^2\mu(B(x,r_k)) \le 8C(1)^{2}\omega'_n r_k^n = \frac{V}{k}}$. Our Inequality~(\ref{bornepreuve}) now reads as
$$	\forall k \ge k_0 \qquad \l_k \le \frac{4C(1)}{r_k^2}= 4C(1)\Big(8 C(1)^{2} \omega'_n \Big)^{2/n}\left(\frac{k}{V}\right)^{2/n}\ .$$
Now if $k < k_0$, then we obviously have $\l_k \le \l_{k_0}$, so that we straightly obtain
\begin{equation}\label{asympt}
 \forall k \in \N^* \qquad \l_k \le \l_{k_0}  + B_n \left(\frac{k}{V}\right)^{2/n} \ ,
 \end{equation}
where we have set $B_n := 4C(1)\Big(8 C(1)^{2}\omega'_n \Big)^{2/n}$. The last thing to do is to estimate the particular eigenvalue $\l_{k_0}$.
\begin{itemize}
	\item[1)] If $k_0=1$, then $\l_{k_0}=\l_1=0$ and we get Inequality~(\eq{riccineg}), with $A_n=0$.\\
	\item[2)] On the contrary, if $k_0\ge 2$, then we deduce $\frac{V}{8C(1)^2\omega'_n } < k_0 \le 2 \frac{V}{8C(1)^2\omega'_n }$. We can apply Inequality~(\eq{asympt}) with $k=k_0$, which implies
	$$	\l_{k_0} \le \frac{4C(1)}{r_{k_0}^2}= 4C(1) 2^{2/n}\ ,$$
	and then Inequality~(\eq{asympt}) is nothing but Inequality~(\eq{riccineg}) with $A_n=4C(1) 2^{2/n}$, $B_n = 
	4C(1)\Big(8 C(1)^{2}\omega'_n \Big)^{2/n}$ and  $a=1$.
	\end{itemize} 
 \qed

\begin{Rq}
For the case $a=0$, a slightly better constant $B_n$ can be otained by making a direct proof instead of plugging $a=0$ in Inequality~(\eq{riccineg}).
\end{Rq}

\medskip
\noindent
\textsc{\textbf{Acknowledgements.}} We thank G.~Carron for helpful comments. DM is grateful to the Mathematical Institute of Neuchâtel for financial support during the redaction of a large part of this text.


\end{document}